\documentclass[letterpaper, 10 pt, conference]{ieeeconf}  

\IEEEoverridecommandlockouts                              
\usepackage{amsthm}
\usepackage{amsmath, amsfonts, amssymb}
\usepackage{mathrsfs, float}
\usepackage{algorithm}
\usepackage{algorithmic}
\usepackage{comment}
\usepackage{bm}
\usepackage{verbatim}
\usepackage{mathtools}
\usepackage{bbm}
\usepackage{dsfont}
\usepackage{booktabs}
\usepackage{relsize}
\usepackage{caption}
\usepackage{subcaption}
\usepackage[english]{babel}
\usepackage{xcolor}
\usepackage{cite}
\usepackage{caption}

\let\labelindent\relax
\usepackage{enumitem}

\makeatletter
\let\NAT@parse\undefined
\makeatother
\usepackage[colorlinks]{hyperref}
\hypersetup{
	colorlinks   = true, 
	urlcolor     = blue, 
	linkcolor    = red, 
	citecolor    = green!50!black   
}

\allowdisplaybreaks

\theoremstyle{plain}

\newtheorem{assumption}{Assumption}
\newtheorem{lemma}{Lemma}

\newtheorem{theorem}{Theorem}

\DeclareMathOperator*{\argmin}{arg\,min}

\def\bA{\bm{A}}
\def\by{\bm{y}}

\def\bw{\bm{w}}

\def\bg{\bm{g}}
\def\bG{\bm{G}}

\def\err{{\cal E}}
\def\errMax{{\cal E}^{\max}}

\def\bxi{\bm{\xi}}

\def\E{\mathbb{E}}

\def\ConstLemOne{{\rm C}_L}

\title{\LARGE \bf
	Linear Speedup of Incremental Aggregated Gradient Methods on Streaming Data
}

\author{Xiaolu Wang$^{1}$, Cheng Jin$^{2}$, Hoi-To Wai$^{1}$, Yuantao Gu$^{2}$
	\thanks{$^{1}$Xiaolu Wang and Hoi-To Wai are with System Engineering \& Engineering Management, Faculty of Engineering, The Chinese University of Hong Kong, Hong Kong SAR. Emails:~\url{xwangcu@gmail.com}, \url{htwai@se.cuhk.edu.hk}.
		$^{2}$Cheng Jin and Yuantao Gu are with Department of Electronic Engineering, Tsinghua University, Beijing. Emails:~\url{jinc21@mails.tsinghua.edu.cn}, \url{gyt@tsinghua.edu.cn}.
		This work is partly supported by CUHK Direct Grant \#4055208.}
}

\begin{document}
	
	\maketitle
	\thispagestyle{empty}
	\pagestyle{empty}
	
	\begin{abstract}
		This paper considers a type of incremental aggregated gradient (IAG) method for large-scale distributed optimization. The IAG method is well suited for the parameter server architecture as the latter can easily aggregate potentially staled gradients contributed by workers. Although the convergence of IAG in the case of deterministic gradient is well known, there are only a few results for the case of its stochastic variant based on streaming data. Considering strongly convex optimization, this paper shows that the streaming IAG method achieves linear speedup when the workers are updating frequently enough, even if the data sample distribution across workers are heterogeneous. We show that the expected squared distance to optimal solution decays at ${\cal O}( (1+T) / (nt) )$, where $n$ is the number of workers, $t$ is the iteration number, and $T/n$ is the update frequency of workers. {Our analysis involves careful treatments of the conditional expectations with staled gradients and a recursive system with both delayed and noise terms, which are new to the analysis of IAG-type algorithms.} Numerical results are presented to verify our findings. 
	\end{abstract}
	
	\section{INTRODUCTION}
	Distributed optimization is an important algorithmic paradigm that has received immense attention due to its wide applicability in machine learning, signal processing, control, etc \cite{yang2019survey,chang2020distributed,li2020federated}. It is suitable for a broad range of circumstances where data are dispersed across multiple entities, e.g., CPU cores, computing clusters, wireless sensors, and wearable devices \cite{yang2019survey}. 
	Classical distributed optimization deals with the case when each worker holds a fixed set of local data samples that is available at any time, which is also referred to as the \emph{batch data learning} setting. 
	However, with the growing scenarios including federated learning \cite{li2020federated} where the data are acquired in an online fashion (e.g., online review and social network platforms) and each data sample is allowed to be used only once \cite{chang2020distributed}, it is important to adapt the distributed optimization algorithms to the streaming data setting.
	
	This paper is concerned with the following stochastic optimization problem:
	\begin{equation}\label{eq:opt}
		\begin{split}
			\min_{\bw\in\mathbb{R}^{d}}~
			&  \frac{1}{n} \sum_{i=1}^{n} F_i(\bw),~~ F_i(\bw) := \mathbb{E}_{ \bxi_i \sim \mathcal{D}_i } \left[ f_i (\bw;\bxi_i) \right],
		\end{split}
	\end{equation}
	where ${\cal D}_i$ represents the data distribution supported on the sample space $\Xi_i$ accessible by worker $i$. 
	The optimization problem shall be solved cooperatively by $n$ workers. For $i \in [n]$, $\bxi_i \in \Xi_i$, the sampled local loss function $f_i( \cdot ; \bxi_i )$ is continuously differentiable and is known to the $i$th worker. We denote by $F( \bw) := (1/n) \sum_{i=1}^n F_i(\bw )$ the global objective function and assume that $F( \bw )$ is strongly convex. 
	
	Consider solving \eqref{eq:opt} in a distributed fashion under the coordination of a central server communicating with the $n$ workers. Each of the worker has access to an independent streaming data source ${\cal D}_i$ and the latest iterate  $\bw$ stored at the server, with which it computes stochastic estimates of the gradient $\nabla F_i( \bw )$. 
	We concentrate on an asynchronous setting where workers can be idle in some iterations, due to, e.g., network connection failure. Notice that using direct average of the stochastic gradients may result in a non-converging algorithm unless the local loss satisfies some form of similarity conditions; see the related studies on FedAvg in \cite{mcmahan2017communication, karimireddy2020scaffold}. Remedies such as designing stochastic control variate have been proposed, e.g., \cite{karimireddy2020scaffold, mishchenko2022proxskip}. 
	
	To deal with the worker asynchrony issue over heterogeneous data, this paper utilizes a parameter server (PS) architecture \cite{assran2020advances, aytekin2016analysis} for distributed optimization of \eqref{eq:opt}. In this setting, the PS maintains a buffer that stores gradient information reported from the workers. The coordinating server then aggregates the information stored in the buffer to update its iterate. While performing gradient aggregation seems to alleviate the reliance on data similarity, we note the aggregated information may contain staled gradients due to worker asynchrony that can affect convergence. 
	To this end, this work inquires the following questions: 
	
	\emph{Does the above distributed algorithm solve \eqref{eq:opt}? Can it achieve linear speedup in convergence rate compared to {centralized/sequential} SGD solving \eqref{eq:opt} taking one sample per iteration?} 
	
	\begin{figure}[t]
		\centering
		\includegraphics[width=0.8\columnwidth]{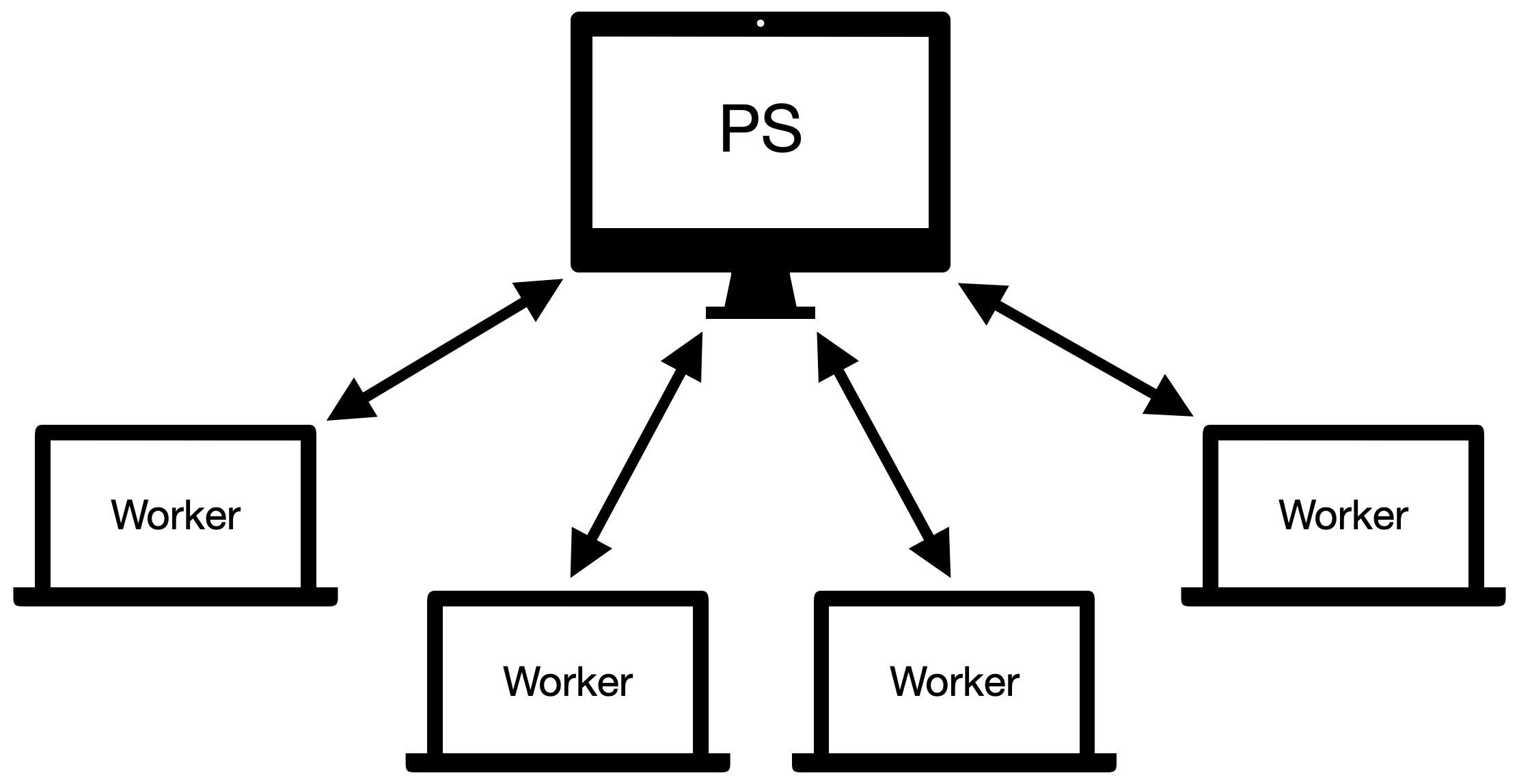}
		\caption{PS architecture. In this architecture, the PS keeps a buffer of size $n \times d$ that stores the latest gradient computed and sent by the workers and distributes the latest $\bw^t$ to the workers.}
		\label{fig:ps}\vspace{-.3cm}
	\end{figure}
	
	\vspace{.1cm}
	\noindent \textbf{Our Contributions.}
	In this paper, we provide an affirmative answer to the above questions. We study a streaming incremental aggregated gradient (sIAG) method adapted from the incremental aggregated gradient (IAG) method \cite{blatt2007convergent} to handle streaming data. Our key results are as follows:
	\begin{itemize}[leftmargin=*]
		\item We show that the sIAG method converges in expectation to the global optimum solution of \eqref{eq:opt} under mild conditions and \emph{without requiring an explicit similarity condition} between $F_i(\bw)$.
		\item Suppose that the maximum stateness of the aggregated gradient is $T$, we show that the expected squared error between the $t$th iterate and optimal solution to \eqref{eq:opt} is ${\cal O}( \frac{1}{t} \frac{ \sigma^2 (1+T) }{n} )$ for the sIAG method, where $\sigma^2$ is the variance of a gradient sample. As such, the sIAG method achieves linear speedup despite that it utilizes staled gradient information with asynchronous workers.
	\end{itemize}
	We remark that our analysis utilizes new analysis techniques to yield the tight bound for sIAG; see Sec.~\ref{sec:ana}. 
	
	\vspace{.1cm}
	\noindent \textbf{Related Works.}
	The studies of IAG method and its variants based on batch data have generated substantial interests after  \cite{blatt2007convergent}. Notably, \cite{feyzmahdavian2014delayed, gurbuzbalaban2017convergence,vanli2018global} have worked on analyzing the convergence rate of IAG under the bounded delay assumption. Extensions have been considered to speed up convergence, e.g., \cite{chen2018lag, wu2022delay} studied adaptive strategies for gradient aggregation, \cite{wai2020accelerating} utilized local Hessian information. 
	
	On the other hand, the study of IAG-like methods under the streaming data setting has received attention from the machine learning community. The closest work to ours is \cite{lian2015asynchronous} that studied an sIAG method with smooth objective function under stronger assumption than ours, e.g., drawing independent samples for staled gradient. 
	As mentioned, federated learning algorithms such as FedAvg \cite{mcmahan2017communication}, FedProx \cite{li2020federated-a}, SCAFFOLD \cite{karimireddy2020scaffold} adopted similar aggregation technique for the local information reported from the workers.

	\section{STREAMING IAG METHOD}
	This section introduces the sIAG method and discusses its implementation in a distributed optimization setting with the PS. To motivate the algorithmic idea of sIAG method, below we first briefly review the IAG method for \eqref{eq:opt} with batch data \cite{blatt2007convergent}. 
	
	Consider a distributed computing architecture with a PS and $n$ workers \cite{assran2020advances}; 
	see Fig.~\ref{fig:ps}. The IAG method assumes that each worker $i \in [n]$ has \emph{full access} to an oracle that queries the gradient of its local loss function $F_i( \bw )$ at any point $\bw \in \mathbb{R}^d$, e.g., when the local dataset is fixed.
	The workers send their computed gradient to the PS that keeps a buffer of $n$ gradient vectors storing the most recent copy of the computed gradient from each worker. Note that during most iterations, this buffer may contain \emph{staled gradient} when the workers is idle at the current iteration. 
	
	To fix notations, we let ${\cal A}_t \subseteq [n]$ to be set of active workers at iteration $t$ and define 
	\begin{equation} \label{eq:tau-def}
		\tau_i(t) = \max\big\{ \tau : \tau \leq t,~i \in {\cal A}_\tau \big\}.
	\end{equation}
	In other words, $\tau_i(t) \leq t$ indicates the iteration number in which the gradient stored in the PS from worker $i$ is computed. 
	In iteration $t \geq 0$, the coordinating server performs the update: 
	\begin{align}
		\bw^{t+1} = \bw^{t} - \frac{\eta}{n} \sum_{i=1}^{n} \nabla F_i (\bw^{\tau_i(t)}) ,
		\label{eq:iag}
	\end{align}
	where $\eta > 0$ is the step size. 
	Observe that the server directly aggregates all local (possibly staled) gradients in the PS as the descent direction. 
	
	Despite that \eqref{eq:iag} utilizes some staled gradients in the updates, a key result established in \cite{gurbuzbalaban2017convergence, vanli2018global} is that IAG admits \emph{linear convergence} towards the optimal solution of \eqref{eq:opt} when $F(\bw)$ is strongly convex and smooth, similar to a \emph{centralized} gradient method for \eqref{eq:opt} under the same setting. Especially, this convergence rate holds regardless of the differences between the local loss functions $F_i(\bw)$. In comparison, without the aggregation step performed with the PS, the distributed algorithm that relies on aggregating only the current gradients reported by ${\cal A}_t$ may converge sublinearly and requires further modification such as taking a diminishing step size \cite{bertsekas2011incremental}.\vspace{.2cm}
	
	\algsetup{indent=1em}
	\begin{algorithm}[t]
		\caption{sIAG Algorithm} 
		\label{alg:siag}
		\begin{algorithmic}[1]
			\STATE {\bf input:} Initialization $\bw^0$, step sizes $\{ \eta_t\}_{t\geq 0}$. 
			\STATE At the PS, initialize the buffer with $\bg_i^{-1} = {\bm 0}$ for $i=1,\ldots,n$. 
			\FOR{$t=0,1,2,\ldots$}
			\STATE A set of workers ${\cal A}_t \subseteq [n]$ is selected/active. 
			\FOR{each worker $i \in {\cal A}_t $}
			\STATE Take $\bw^t$ from the PS and draw a sample  $\bxi_i^t \sim {\cal D}_i$.
			\STATE Compute the stochastic gradient $\nabla f_i ( \bw^t; \bxi_i^t )$ and send back to the PS.
			\ENDFOR
			\STATE Update the buffer as 
			\begin{align*}
				&\bg_i^t \leftarrow \nabla f_i ( \bw^t; \bxi_i^t )~\text{for}~i \in {\cal A}_t,
				\\
				&\bg_i^t \leftarrow \bg_i^{t-1}~\text{for}~i \notin {\cal A}_t \\[-.9cm]
			\end{align*}
			\STATE Compute sIAG update: $\bw^{t+1} \leftarrow \bw^t - (\eta_t/n) \sum_{i=1}^n \bg_i^t$.
			\ENDFOR
		\end{algorithmic}
	\end{algorithm}
	
	\noindent \textbf{Streaming IAG Method.} We are interested in a variant of the IAG method utilizing \emph{streaming data}. We consider the generic form for \eqref{eq:opt} where each of the local loss function is (possibly) stochastic. Unlike the IAG method, the streaming IAG (sIAG) method considers that each worker only has access to a \emph{stochastic oracle} that queries an unbiased and independent estimate for the gradient of local loss function, denoted by 
	\begin{equation}
		\nabla f_i( \bw; \bxi_i ),~~\bxi_i \sim {\cal D}_i.
	\end{equation}
	Note that $\E[ \nabla f_i( \bw; \bxi_i ) ] =  \nabla F_i( \bw )$.
	The sIAG method reads 
	\begin{align}
		\bw^{t+1} = \bw^{t} - (\eta_t/n) \, \bg^t,
		\label{eq:siag}
	\end{align}
	for any $t \geq 0$, where $\eta_t > 0$ is a (possibly time varying) step size and  
	\begin{align}
		& \bg^t	\coloneqq \sum_{i=1}^{n} \nabla f_i ( \bw^{\tau_{i}(t)}; \bxi_i^{\tau_{i}(t)} )
		\label{eq:gt}
	\end{align}
	denotes the aggregated \emph{stochastic gradients}. Note that in the above, we adopted the same notations as in the IAG method where the index $\tau_i(t)$ was defined in \eqref{eq:tau-def}.
	
	The sIAG method \eqref{eq:siag}, \eqref{eq:gt} is motivated by an instantaneous gradient computation model that can be readily implemented on the parameter server architecture; as summarized in Algorithm~\ref{alg:siag}. 
	Observe that the algorithm requires active workers to return the stochastic gradient before the current iteration concludes. 
	This restriction, while mild as only stochastic gradients are required, might be relaxed by allowing further delays between the sampled gradient and the iterate in which it is computed. However, in the interest of space, we focus on a current simplified setting.
	
	We recall that as $\tau_i(t) = t$ for $i\in\mathcal{A}_t$ and $\tau_i(t) < t$ for $i\notin\mathcal{A}_t$, we may express the aggregated gradient \eqref{eq:gt} alternatively as an incremental update via
	\begin{align*}
		&\bg^t 
		= \bg^{t-1} \hspace{-1mm}-\hspace{-1.5mm} \sum_{i\in\mathcal{A}_t} \hspace{-1mm}\nabla f_i ( \bw^{\tau_{i}(t-1)}; \bxi_i^{\tau_{i}(t-1)} ) \hspace{-1mm}+\hspace{-1.5mm} \sum_{i\in\mathcal{A}_t}\hspace{-1mm} \nabla f_i ( \bw^t; \bxi_i^t ) .
	\end{align*}
	We remark that the above update recursion for $\bg^t$ is related to the famous SAG method \cite{schmidt2017minimizing}. It is reduced to the SAG method when $f_i \equiv f$, $\bxi_i^t = \bxi^t$, and ${\cal D}_i \equiv {\cal D}$ has a finite support for all $i \in [n]$. 
	
	While convergence guarantees for the IAG method \eqref{eq:iag} are well known, the stochastic variant, sIAG, considered in \eqref{eq:siag}, \eqref{eq:gt} has received less attention. As the first step towards understanding the behavior of sIAG, the next section analyzes the convergence rate of sIAG under the standard setting when $F(\bw)$ is strongly convex and smooth.

	\section{CONVERGENCE ANALYSIS} \label{sec:ana}
	We shall begin by stating some assumptions on \eqref{eq:opt} and the sIAG method that are necessary for our analysis. First,	
	\begin{assumption}\label{as:function}
		Problem \eqref{eq:opt} satisfies: {\sf (i)} $F(\bw)$ is $\mu$-strongly convex with $\mu > 0$; {\sf (ii)} for any $i \in [n]$, the gradient $\nabla f_i( \bw )$ is $L$-Lipschitz continuous.
	\end{assumption}
	\noindent The above specifies the function class of interest for \eqref{eq:opt}. We define $\bw^\star = \argmin_{ \bw \in \mathbb{R}^d} F(\bw)$. We also assume:
	\begin{assumption}\label{as:sigma}
		There exists a constant $\sigma \geq 0$ such that
		\begin{equation}
			\E \| \nabla f_i(\bw;\bxi_i) - \nabla F_i(\bw) \|_2^2 \leq \sigma^2 ( 1 + \| \bw - \bw^\star \|^2 ),
		\end{equation}
		for all $i \in [n]$, $\bw \in \mathbb{R}^d$.
	\end{assumption}
	\begin{assumption}\label{as:delay}
		There exists $T \geq 0$ such that 
		\begin{equation}
			\tau_i(t) \geq t-T,~\forall~i \in [n], ~t \geq 0.
		\end{equation}
	\end{assumption}
	\noindent The above assumptions state that the stochastic gradients computed in sIAG has bounded variance, and the staled gradient delay is bounded by $T$. Notice that $n / T$ corresponds roughly to the number of workers that are active at any iteration.  
	
	We define respectively the following notations for the $t$th suboptimality gap and its delayed version: 
	\begin{align}
		\err_t \coloneqq \E \|\bw^t - \bw^* \|_2^2, \quad
		\errMax_t \coloneqq \max_{s\in[ (t-2T)_+,t]} \err_s,
	\end{align}
	where $x_+ \coloneqq \max\{x,0\}$ for $x\in\mathbb{R}$.
	Furthermore, it is instrumental to define the following filtration:
	\begin{equation}
		{\cal F}_t \coloneqq \sigma( \bxi_i^s , i \in [n], s = 0,\ldots, t-1 ),
	\end{equation}
	where $\sigma(\cdot)$ denotes the sigma algebra generated by the random variables in the operand. 
	Observe that for any $t \geq 0$, $\bw^t$ is measurable with respect to (w.r.t.) ${\cal F}_t$. We shall use the shorthand notation $\E_t[\cdot] := \E[ \cdot | {\cal F}_t ]$ for conditional expectation.
	
	From \eqref{eq:siag}, we deduce that for any $t \geq 0$,
	\begin{equation} \label{eq:ana-iter}
		\begin{split}
			\err_{t+1} & \hspace{-0.5mm}=\hspace{-0.5mm} \err_t \hspace{-0.5mm}-\hspace{-0.5mm} 2 \eta_t \mathbb{E} \hspace{-0.5mm}\left[\hspace{-0.5mm} \left\langle\hspace{-1mm} \bw^t \hspace{-0.5mm}-\hspace{-0.5mm} \bw^\star , \frac{1}{n} \bg^t \hspace{-0.6mm}\right\rangle \hspace{-0.5mm}\right] \hspace{-1mm}+\hspace{-0.5mm} \eta_t^2 \mathbb{E} \left[ \left\| \frac{1}{n} \bg^t \right\|^2 \right]\hspace{-0.5mm}.
		\end{split}
	\end{equation}
	We shall control the last two terms in the above decomposition. Observe the following lemmas:
	\begin{lemma}\label{lem:aa}
		Suppose that Assumptions \ref{as:function}--\ref{as:delay} hold. Then, 
		\begin{align}
			& \E \| (1/n) \bg^t \|_2^2 \leq \frac{ 2 \sigma^2 } {n} + \ConstLemOne
			\max_{ s \in [ (t-T)_+, t ] } \err_s,
			\label{eq:gt-bound}
		\end{align}
		for any $t \geq 0$,
		where $\ConstLemOne \coloneqq 20 L^2 + \frac{2\sigma^2}{n}$.
	\end{lemma}
	\noindent The above is a natural consequence of \eqref{eq:gt}, where the bound \eqref{eq:gt-bound} depends on both the stochastic gradient variance $\sigma^2$ and the delayed optimality gap $\max_{ s \in [ (t-T)_+, t ] } \err_s \leq \errMax_t$. 
	The detailed proof can be found in Appendix~\ref{ap:aa}.
	
	The next lemma controls the inner product term in \eqref{eq:ana-iter}:
	\begin{lemma}\label{lem:bb}
		Suppose that Assumptions \ref{as:function}--\ref{as:delay} hold. Then
		\begin{align}
			&\E \left\langle \bw^t - \bw^*, (1/n) \, \bg^t \right\rangle
			\label{eq:innerprod-bound}
			\\
			&\geq \frac{\mu}{4} \err_t - \left[ \ConstLemOne T \eta_{t-T} + \left( \frac{\mu}{4} + \frac{5L^2}{2 \mu} \right) \ConstLemOne T^2 \eta_{t-T}^2 \right] \err_{t}^{\max}
			\nonumber
			\\
			&\quad -\left[ 2T\eta_{t-T} + \left( \frac{\mu}{2} + \frac{5L^2}{\mu} \right) T^2 \eta_{t-T}^2 \right] \frac{\sigma^2}{n},~\forall~t \geq 0.
			\nonumber
		\end{align}
	\end{lemma}
	\noindent The above lemma shows that the inner product term is \emph{lower bounded} by $(\mu/4) \err_t$ with (negative) perturbation terms that are controllable by the step sizes. 
	
	Before substituting \eqref{eq:gt-bound}, \eqref{eq:innerprod-bound} into \eqref{eq:ana-iter} to derive the convergence rate of sIAG method, we highlight that Lemma~\ref{lem:bb} deviates from the standard analysis for SGD with strongly convex objective function. In the case of standard SGD where $\bg^t = \nabla f( \bw^t ; \bxi^t )$ with $\E_t[ \bg^t ] = \nabla F( \bw^t )$, the law of total expectation and the strong convexity of $F$ imply that
	\begin{align*}
		&\E \langle \bw^t - \bw^\star , \bg^t \rangle 
		= \E \langle \bw^t - \bw^\star , \E_t [\bg^t] \rangle 
		\nonumber
		\\
		&= \E \langle \bw^t - \bw^\star , \nabla F( \bw^t ) \rangle \geq \mu \err_t.
	\end{align*}
	However, in the sIAG method, $\bg^t$ is a function of $\bxi_i^s$ for $i \in [n]$ and $s=0,1,\dots,t$. Note that for any $i \notin \mathcal{A}_t$, $\bxi_i^{\tau_i(t)}$ is $\mathcal{F}_t$-measurable and thus $\E_t [\nabla f_i(\bw_i^{\tau_i(t)};\bxi_i^{\tau_i(t)})] \neq \nabla F_i(\bw^{\tau_i(t)})$. This implies that $\bg^t$ is \emph{not} independent of $\mathcal{F}_t$ and thus
	\begin{align*}
		\E \left\langle \bw^t - \bw^\star, \bg^t \right\rangle  
		&\neq \E \left\langle \bw^t - \bw^\star , \sum_{i=1}^n \nabla F_i( \bw^{\tau_i(t)} ) \right\rangle.
	\end{align*}
	Our remedy is to consider\footnote{We remark that \cite{lian2018asynchronous} considered an algorithm that involve similar staled aggregation property to sIAG but have employed a simplifying assumption that $\bw^t$ is independent of $\bg^t$. Leveraging this property, their algorithm achieves linear speedup regardless of $T$. We conjecture that such linear speedup cannot be obtained when considering the realistic conditions for sIAG method that $\bw^t$ is not independent of $\bg^t$.} the following decomposition for the inner product:
	\begin{equation} \label{eq:ana-split}
		\begin{split} 
			& \langle \bw^t - \bw^*, \nabla f_i( \bw^{\tau_i(t)}; \bxi_i^{\tau_i(t)} ) \rangle \\
			& = \langle \bw^t - \bw^{\tau_i(t)} + \bw^{\tau_i(t)} - \bw^*, \nabla f_i( \bw^{\tau_i(t)}; \bxi_i^{\tau_i(t)} ) \rangle.
		\end{split}
	\end{equation}
	Observe that $\nabla f_i( \bw^{\tau_i(t)}; \bxi_i^{\tau_i(t)} )$ is independent of $\bw^{\tau_i(t)}$ and thus the simplification $\E_{ \tau_i(t)} [ \langle \bw^{\tau_i(t)} - \bw^\star, \nabla f_i( \bw^{\tau_i(t)}; \bxi_i^{\tau_i(t)} ) \rangle ] = \langle \bw^{\tau_i(t)} - \bw^\star, \nabla F_i( \bw^{\tau_i(t)} ) \rangle$.
	
	Furthermore, we develop the following bound to control the size of difference $\| \bw^t - \bw^{\tau_i(t)} \|_2^2$:
	\begin{lemma}\label{lem:diff}
		Suppose that Assumptions \ref{as:function} and \ref{as:sigma} holds and $\{\eta_t\}_{t\geq0}$ is a monotonically non-increasing sequence. Then, it holds for all $i \in [n]$ and $t \geq 0$ that
		\begin{align}
			\E [ \| \bw^t - \bw^{\tau_{i}(t)} \|_2^2 ]
			\leq T^2 \eta_{t-T}^2 \left( \frac{2 \sigma^2}{n} + \ConstLemOne \errMax_{t} \right).
			\label{eq:wt-diff-bound}
		\end{align}
	\end{lemma}
	\noindent Observe that the bound is composed of the stochastic gradient's variance and the (delayed) optimality gap. Note that the proof is obtained as a consequence of Lemma~\ref{lem:aa}. 
	
	Importantly, \eqref{eq:wt-diff-bound} allows us to handle the term $\langle \bw^t - \bw^{\tau_i(t)}, \nabla f_i( \bw^{\tau_i(t)}; \bxi_i^{\tau_i(t)} ) \rangle$ through applying Cauchy-Schwarz inequality. 
	The detailed proofs of Lemma~\ref{lem:bb}, \ref{lem:diff} are found in Appendix~\ref{ap:bb}, \ref{ap:diff}, respectively.  
	
	Substituting Lemmas~\ref{lem:aa} and \ref{lem:bb} into \eqref{eq:ana-iter} directly yields the following recursive system: for any $t \geq 0$,
	\begin{align} 
		& \err_{t+1} \leq (1 - (\mu/2) \, \eta_t) \, \err_t \label{eq:ana-recur} 
		\\
		& + \left[ 1 + 2 T + \left( \frac{\mu}{2} + \frac{5L^2}{\mu} \right) T^2 \eta_{t-T} \right] \ConstLemOne \eta_{t-T}^2 \err_{t}^{\max}
		\nonumber
		\\
		& + \left[ 1 + 2T + \left( \frac{\mu}{2} + \frac{5L^2}{\mu} \right) T^2 \eta_{t-T} \right] \eta_{t-T}^2 \frac{2 \sigma^2}{n}. 
		\nonumber
	\end{align}
	If we ignore the last term proportional to $\sigma^2/n$, then \eqref{eq:ana-recur} reduces into a contracting recursion with delays that has been studied by \cite{feyzmahdavian2014delayed}. Subsequently, $\err_t$ converges to zero at a linear rate when a constant step size is used. 
	
	The introduction of the noise-related terms in \eqref{eq:ana-recur} has led to a new recursive system with delayed terms that have not been covered in prior works. 
	As an attempt to derive a tight bound, we fix the step size as 
	\begin{equation}
		\eta_t = \beta / (t + \gamma)
	\end{equation}
	for some $\gamma, \beta > 0$ and obtain the following convergence rates for $\err_t$ using induction:
	\begin{theorem}\label{thm}
		Suppose that $\beta > 2/\mu$ and 
		\begin{align*}
			\gamma \geq 2T + \max \left\{ \frac{16C_L\beta^2}{\mu\beta-2}\bar{\rho}(T), \sqrt{\frac{8C_L\beta^2}{\mu\beta-4}\bar{\rho}(T)} \right\},
		\end{align*}
		where $\bar{\rho}(T) \coloneqq 1 + 2 T + \left( \frac{\mu}{2} + \frac{5L^2}{\mu} \right) \beta T$.
		Let
		\begin{align*}
			&
			\delta_1 \coloneqq \frac{32 \beta^2 \bar{\rho}(T)}{ \mu \beta - 2 } + 1, ~\delta_2 = \gamma^2 \err_0.
		\end{align*}
		Then, it holds for all $t \geq 0$ that 
		\begin{align}
			&\err_{t} \leq \frac{\delta_1}{\gamma+t} \frac{\sigma^2}{n} + \frac{\delta_2}{(\gamma+t)^2}.
			\label{eq:rate}
		\end{align}
	\end{theorem}
	\noindent The detailed proof can be found in Appendix~\ref{ap:thm}. Theorem \ref{thm} indicates that the expected squared distance to optimal solution decays at ${\cal O}( (1+T) / (nt) )$.
	We note that in the interest of space, the constants in the bound are not fully optimized. In general, obtaining tight bounds for recursive system of the form \eqref{eq:ana-recur} is an interesting open problem of independent interest. 
	
	The bound \eqref{eq:rate} shows that the sIAG method converges (in expectation) towards the optimal solution of \eqref{eq:opt} at the rate of ${\cal O}( \frac{1}{t} \cdot \frac{( 1+T ) \sigma^2 }{n} )$. 
	We notice that: {\sf (A)} the sublinear rate of ${\cal O}(1/t)$ is similar to that of existing analysis with stochastic gradient methods \cite{moulines2011non} and is in the same order of the minimax lower bound \cite{agarwal2012information}, {\sf (B)} the constant factor ${( 1+T ) \sigma^2 } / {n}$ indicates that \emph{linear speedup} can be achieved when the delay satisfies $T = {\cal O}(1)$. This rate is reasonable since the number of samples taken per iteration is approximately $n/T$, the linear speedup ratio should be of the same order. 
	We note that similar slow down due to the delay $T$ is also reported in the analysis for IAG \cite{gurbuzbalaban2017convergence, vanli2018global}.

	
	\section{NUMERICAL SIMULATIONS}
	We evaluate the empirical performance of the sIAG algorithm on synthetic data. We independently generate $n$ parameters $\bw_1^*,\bw_2^*,\dots,\bw_n^*$ according to the uniform distribution on $[0,1]^d$.
	Each data point sampled by worker $i$ takes the form $\bxi_i=(\bA_i,\by_i)$, where $\bA_i\in\mathbb{R}^{p\times d}$ and $\by_i\in\mathbb{R}^p$. The entries of $\bA_i$ are independent and follow the Gaussian distribution ${\cal N}(0,1)$ and $ \by_i \sim \mathcal{N}(\bA_i\bw_i^*,\sigma^2\bm{I}_p)$.
	The loss functions are defined as $f_i(\bw;\bA_i,\by_i)= \frac{1}{2} \lVert\bA_i\bw-\by_i\rVert_2^2$ for $i\in[n]$. If $[ \bA_1; \cdots; \bA_n ]$ is full-rank (which holds almost surely when $nd \geq p$), it is obvious that $\bw^*=\sum_{i=1}^n \bw_i^*/{n}$.
	
	We compare sIAG with \emph{non-aggregated SGD}, which uses the following descent direction at the $t$th iteration: 
	\begin{align} \textstyle
		\bg_{\sf SGD}^t \coloneqq \sum_{ i: \tau_i(t) = t } \nabla f ( \bw^{\tau_{i}(t)}; \bxi_i^{\tau_{i}(t)} ).
		\label{eq:sgd}
	\end{align}
	We simulate three types of worker selection schemes. The first selection scheme chooses one workers at each iteration cyclically, i.e., at iteration $t$, the $(t \, {\rm mod} \, n + 1)$th worker is active. Notice that $T=n$ in this case and there is no linear speedup according to Theorem~\ref{thm}. The second selection scheme chooses the workers \emph{uniformly at random}. It models the scenario when the workers are equally efficient. Further, the workers will be selected at least once in no more than 15 iterations. 
	The third selection scheme adopts \emph{non-uniform selection}. It models a more realistic scenario when the workers are heterogeneous in terms of efficiency. Specifically, worker $i$ is selected at least once in no more than $T_i$ iterations, where $T_i$ is uniformly distributed on $\{10, \ldots, 20\}$. The faster workers are selected more frequently. Notice that $T=15$ in the second and third scheme. 
	
\begin{figure}[t]
		\centering
		\includegraphics[width=0.9\columnwidth]{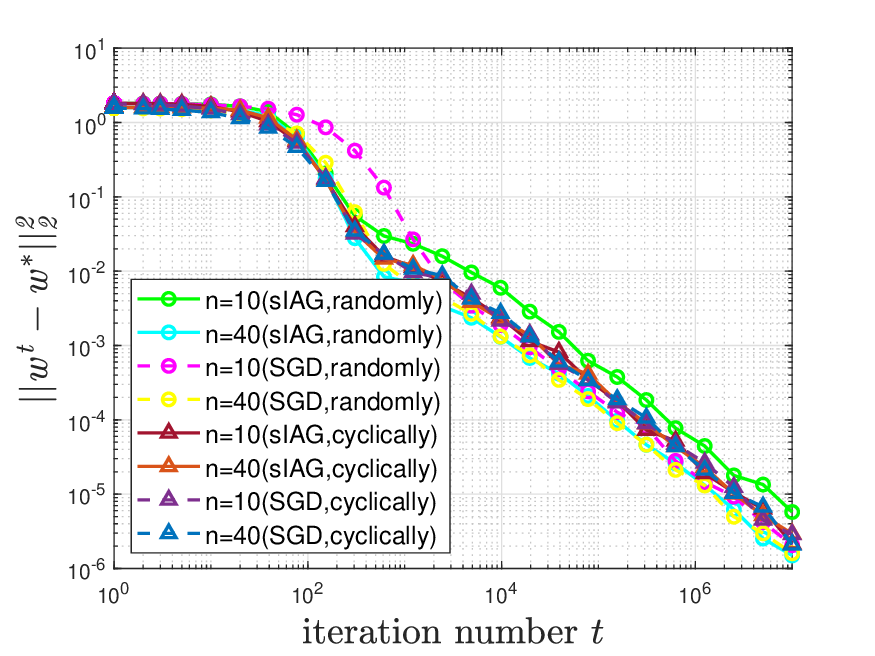}
		\includegraphics[width=0.9\columnwidth]{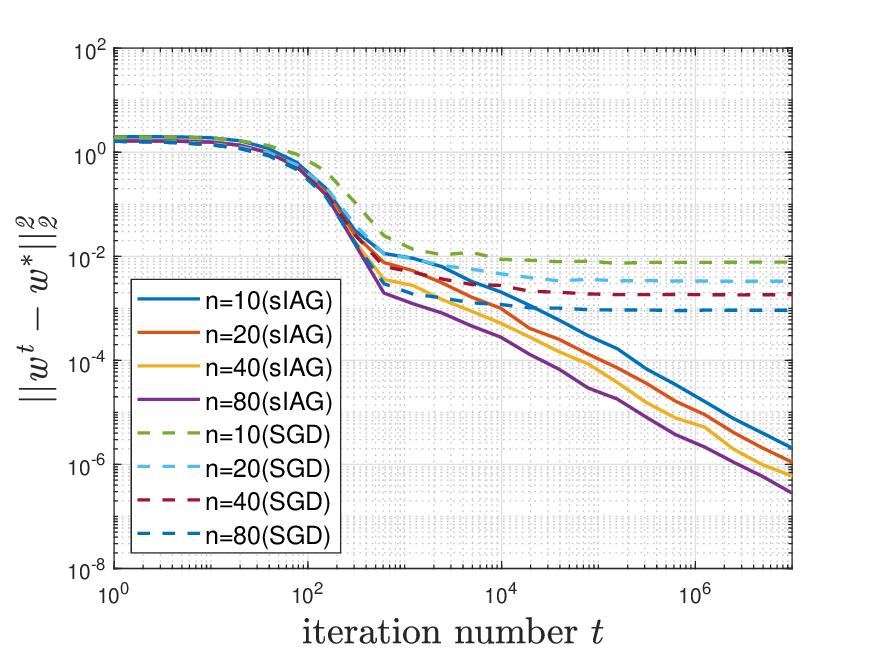}
		\caption{Convergence of sIAG and SGD with: (Top) uniform worker selection; (Bottom) non-uniform worker selection.}\vspace{-.2cm}
		\label{fig:random-uniform}
	\end{figure}

	Fig. \ref{fig:random-uniform} presents the convergence of two algorithms with $d=20$, $p=10$, and $\sigma=0.1$ under the uniform and non-uniform selection schemes. Observe that under the uniform selection scheme, the sIAG and SGD achieve comparable convergence performance and both exhibit linear speedup as the number of agents increases. 
	We also observe that there is no linear speedup with the cyclical selection scheme. 
	On the other hand, under the non-uniform worker selection scheme, the sIAG still enjoys linear speedup while SGD is not converging to an optimal solution, as the non-uniform selection scheme has led to biased stochastic gradient. This validates the necessity of aggregating (possibly) staled gradients in \eqref{eq:siag} as opposed to using only the latest gradients in \eqref{eq:sgd} in the presence of system heterogeneity.
	
	\section{CONCLUSION}
	We proposed the sIAG algorithm for distributed optimization over the parameter server architecture with heterogenous streaming data. The sIAG method is adapted from the classical IAG method on batch data. We established that sIAG achieves linear speedup compared to the sequential SGD for strongly convex problems. Our analysis relies on careful treatments of the conditional expectations with staled gradients (see \eqref{eq:ana-split}) and a new recursive system with both delayed and noise-related terms (see \eqref{eq:ana-recur}), which can be of independent interest. Numerical results on synthetic data verify our theoretical findings and show significant advantages of sIAG over the non-aggregated SGD method when the workers are not uniformly selected.
	
	\bibliographystyle{ieeetr}
	\bibliography{sIAG}
	
	\appendix
	\section{Missing Proofs}
	To simplify notations, throughout this appendix, we denote 
	$\bm{g}_i^{t} \coloneqq \nabla f_i(\bw^{t}; \bxi_i^{t})$
	and $\bG_i^{t} \coloneqq \nabla F_i(\bw^{t})$ as respectively the stochastic gradient and exact gradient of the local loss function for any $i \in [n]$ and $t \geq 0$.

	\subsection{Proof of Lemma \ref{lem:aa}} \label{ap:aa}
	Observe that 
	\begin{align}
		&\E \|(1/n)\bg^t\|_2^2 = \E \left\| \frac{1}{n} \sum_{i=1}^{n} \bm{g}_i^{\tau_{i}(t)} \right\|_2^2
		\label{eq:eee}
		\\
		=& \E \left\| \frac{1}{n} \sum_{i=1}^{n} \left( \bm{g}_i^{\tau_{i}(t)} - \bG_i^{\tau_{i}(t)}	+ \bG_i^{\tau_{i}(t)} \right) \right\|_2^2
		\nonumber
		\\
		\leq& 2 \E \left\| \frac{1}{n} \sum_{i=1}^{n} \left( \bm{g}_i^{\tau_{i}(t)} \hspace{-0.5mm}-\textbf{} \bG_i^{\tau_{i}(t)} \right) \right\|_2^2	\hspace{-2mm}+\hspace{-0.5mm} 2 \E \left\| \frac{1}{n} \sum_{i=1}^{n} \bG_i^{\tau_{i}(t)} \right\|_2^2.
		\notag
	\end{align}
	
	We then upper bound the first term in \eqref{eq:eee}. Assume without loss of generality that for any $i,j \in [n]$ such that $i \neq j$, we have $\tau_{i}(t) < \tau_{j}(t)$. 
	Then, we have $\mathcal{F}_{\tau_i(t)} \subseteq \mathcal{F}_{\tau_j(t)}$.
	Hence, using the law of total expectation and the fact that 
	$\E_{ \tau_j(t) }  \left[ \bm{g}_j^{\tau_j(t)} \right] = \bG_j^{\tau_j(t)}$,
	we obtain
	\begin{align}
		&\E \left\langle \bm{g}_i^{\tau_i(t)} - \bG_i^{\tau_i(t)},  \bm{g}_j^{\tau_j(t)} - \bG_j^{\tau_j(t)} \right\rangle
		\label{eq:0} \\
		& = \E \left[ \left\langle \bm{g}_i^{\tau_i(t)} - \bG_i^{\tau_i(t)},  \E_{\tau_j(t)} \left[ \bm{g}_j^{\tau_j(t)} - \bG_j^{\tau_j(t)} \right] \right\rangle \right] = 0
		\nonumber.
	\end{align}
	Then, it follows from \eqref{eq:0} and Assumptions~\ref{as:sigma}, \ref{as:delay} that
	\begin{align}
		& \E \left\| \frac{1}{n} \sum_{i=1}^{n} ( \bg_i^{\tau_i(t)} - \bG_i^{\tau_i(t)} ) \right\|_2^2 
		\label{eq:ggg} \\
		& = \frac{1}{n^2} \sum_{i=1}^{n} \E \left[ \| \bg_i^{\tau_i(t)} - \bG_i^{\tau_i(t)} \|_2^2 \right]
		\nonumber \\
		& \leq \frac{1}{n^2} \sum_{i=1}^n \sigma^2 ( 1 + \err_{\tau_i(t)} )  
		\leq \frac{\sigma^2}{n} \left( 1 + \max_{s \in [(t-T)_+, t]} \err_s \right). \nonumber 
	\end{align}
	Besides, let $\bG_i^\star \coloneqq \nabla F_i(\bw^*)$ and observe that $\sum_{i=1}^n \bG_i^\star = {\bm 0}$, then we upper bound the second term in \eqref{eq:eee} as follows:
	\begin{align}
		&\E \left\| \frac{1}{n} \sum_{i=1}^{n} \bG_i^{\tau_{i}(t)} \right\|_2^2 
		= \E \left\| \frac{1}{n} \sum_{i=1}^{n} \left( \bG_i^{\tau_{i}(t)} - \bG_i^t + \bG_i^t \right) \right\|_2^2
		\nonumber\\
		&\leq \frac{2}{n^2} \E \left\| \sum_{i=1}^{n} \left( \bG_i^{\tau_{i}(t)} - \bG_i^t \right) \right\|_2^2 + \frac{2}{n^2} \E \left\| \sum_{i=1}^{n} \bG_i^t \right\|_2^2
		\nonumber\\
		&\leq \frac{2 L^2}{n} \sum_{i=1}^{n} \E \| \bw^t - \bw^{\tau_{i}(t)} \|_2^2 + \frac{2 L^2}{n} \sum_{i=1}^{n} \E \| \bw^t - \bw^\star \|_2^2
		\nonumber\\
		&=  6L^2 \E \| \bw^t - \bw^* \|_2^2 + \frac{4 L^2}{n} \sum_{i=1}^{n} \E \| \bw^{\tau_i(t)} - \bw^\star \|_2^2 \nonumber \\ 
		&\leq  10L^2 \max_{ s \in [ (t-T)_+, t ] } \err_s,
		\label{eq:ggg2}
	\end{align}
	Plugging \eqref{eq:ggg} and \eqref{eq:ggg2} back into \eqref{eq:eee} gives the desired bound.
	
	\subsection{Proof of Lemma \ref{lem:bb}} \label{ap:bb}
	Let $\bar{\tau}(t) \coloneqq \min_{i\in[n]} \tau_i(t)$, then we have
	\begin{align}
		&\E \left\langle \bw^t - \bw^*, (1/n) \bg^t \right\rangle
		\nonumber
		\\
		=& \underbrace{\E\left\langle \bw^t -\bw^{\bar{\tau}(t)}, \frac{1}{n} \bg^t \right\rangle}_{A}+ \underbrace{\E\left\langle\bw^{\bar{\tau}(t)}-\bw^*, \frac{1}{n} \bg^t \right\rangle}_{B}.
		\label{eq:club+spade}
	\end{align}
	Note that as explained previously, $\bg^t$ is independent of $\bw^{\bar{\tau}(t)}$ conditioning on $\mathcal{F}_{\bar{\tau}(t)}$ and thus the $B$ term can be controlled easily. 
	In the sequel, we bound the $A,B$ terms in \eqref{eq:club+spade} to obtain desired result.
	\vspace{.1cm}
	
	\noindent
	\underline{i) Bounding $A$:} Notice that $\bg^t$ not conditionally \emph{independent} of $\bw^t$. Our remedy is to observe the decomposition $\bw^t - \bw^{\bar{\tau}(t)} = \sum_{s=\bar{\tau}(t)}^{t-1} (\eta_s/n) \bg^s$. We can lower bound $A$ as follows:
	\begin{align}
		&A = \E \left\langle \sum_{s=\bar{\tau}(t)}^{t-1} \eta_s \bg^s, \bg^t \right\rangle = \sum_{s=\bar{\tau}(t)}^{t-1} \eta_s \E \left\langle \bg^s, \bg^t \right\rangle
		\nonumber\\
		&\geq \sum_{s=\bar{\tau}(t)}^{t-1} \eta_s \left( - \frac{1}{2} \E \| \bg^s \|_2^2 - \frac{1}{2} \E \| \bg^t \|_2^2 \right)
		\nonumber
		\\
		&\geq - T \eta_{t-T} \left( \frac{2 \sigma^2}{n} + \ConstLemOne \err_t^{\max} \right),
		\label{eq:club}
	\end{align}
	where the last inequality uses Lemma~\ref{lem:aa} and the fact that $t-\bar{\tau}(t) \leq T$, $\eta_{t-T} \geq \eta_s$ for any $s = \bar{\tau}_t, \ldots, t-1$.\vspace{.1cm}
	
	\noindent
	\underline{ii) Bounding $B$:} Since $\bar{\tau}(t) \leq \tau_i(t)$ for all $i \in [n]$, we have  $\mathcal{F}_{\bar{\tau}(t)} \subseteq \mathcal{F}_{\tau_i(t)}$. Thus, using $\bg^t = \sum_{i=1}^{n} \bg_i^{\tau_i(t)}$ and the law of total expectation, we have
	\begin{align}
		B &= \E \left[ \E_{\bar{\tau}(t)} \left[ \left\langle \bw^{\bar{\tau}(t)} - \bw^*, \frac{1}{n} \sum_{i=1}^{n} \bg_i^{\tau_i(t)} \right\rangle \right] \right]
		\nonumber
		\\
		&= \underbrace{\E \left\langle \bw^{\bar{\tau}(t)} - \bw^t, \frac{1}{n} \sum_{i=1}^{n} \bG_i^{\tau_i(t)} \right\rangle}_{B_1} 
		\nonumber\\
		&\quad + \underbrace{\E \left\langle \bw^t - \bw^*, \frac{1}{n} \sum_{i=1}^{n} \bG_i^{\tau_i(t)} \right\rangle}_{B_2}.
		\label{eq:star+tri}
	\end{align}
	Since $\sum_{i=1}^n \bG_i^* = \bm{0}$, for any  $\alpha > 0$, we have
	\begin{align}
		&B_1 = \E \left\langle \bw^{\bar{\tau}(t)} - \bw^t, \frac{1}{n} \sum_{i=1}^{n} \left( \bm{G}_i^{\tau_i(t)} - \bm{G}_i^t \right) \right\rangle
		\nonumber
		\\
		&\qquad + \E \left\langle \bw^{\bar{\tau}(t)} - \bw^t, \frac{1}{n} \sum_{i=1}^{n} \left( \bm{G}_i^t - \bm{G}_i^* \right) \right\rangle
		\nonumber
		\\
		&\geq - \frac{\alpha}{2} \E \| \bw^{\bar{\tau}(t)} - \bw^t \|_2^2 - \frac{1}{2\alpha} \E \left\| \frac{1}{n} \sum_{i=1}^{n} ( \bm{G}_i^{\tau_i(t)} - \bm{G}_i^t ) \right\|_2^2 
		\nonumber
		\\
		&\quad - \frac{\alpha}{2} \E \| \bw^{\bar{\tau}(t)} - \bw^t \|_2^2 - \frac{1}{2\alpha} \E \left\| \frac{1}{n} \sum_{i=1}^{n} \left( \bm{G}_i^t - \bm{G}_i^* \right) \right\|_2^2
		\nonumber
		\\
		&\geq -\alpha \E \| \bw^{\bar{\tau}(t)} - \bw^t \|_2^2
		\nonumber
		\\
		&\quad - \frac{L^2}{2n\alpha} \sum_{i=1}^{n} \E \| \bw^{\tau_{i}(t)} - \bw^t \|_2^2 - \frac{L^2}{2n\alpha} \sum_{i=1}^{n} \E \| \bw^t - \bw^* \|_2^2.
		\nonumber
	\end{align} 
	Invoking Lemma~\ref{lem:diff} allows us to further lower bound the above terms by:
	\begin{align} 
		B_1 &\geq -\alpha T^2 \eta_{t-T}^2 \left( \frac{2 \sigma^2}{n} + \ConstLemOne \errMax_{t} \right)
		\nonumber
		\\
		&\quad -\frac{L^2 T^2}{2 \alpha} \eta_{t-T}^2 \left( \frac{2 \sigma^2}{n} + \ConstLemOne \errMax_t \right)  - \frac{L^2}{2\alpha} \err_t
		\nonumber
		\\
		&= -\left( 2\alpha + \frac{L^2}{\alpha} \right) T^2 \eta_{t-T}^2 \frac{\sigma^2}{n} 
		\nonumber
		\\
		&\quad - \left( \alpha + \frac{L^2}{2 \alpha} \right) \ConstLemOne T^2 \eta_{t-T}^2 \err_{t}^{\max} - \frac{L^2}{2\alpha} \err_t,
		\label{eq:bigstar}
	\end{align}
	On the other hand, the term $B_2$ can be controlled as
	\begin{align}
		B_2 &= \E \left\langle \bw^t - \bw^*, \frac{1}{n} \sum_{i=1}^{n} ( \bG_i^{\tau_{i}(t)} - \bG_i^t ) \right\rangle 
		\nonumber
		\\
		&\qquad + \E \left\langle \bw^t - \bw^*, \frac{1}{n} \sum_{i=1}^{n} \bG_i^t \right\rangle
		\nonumber
		\\
		&\geq - \frac{\mu}{2} \E \| \bw^t - \bw^* \|_2^2 - \frac{1}{2\mu} \E \left\| \frac{1}{n} \sum_{i=1}^{n} ( \bG_i^{\tau_{i}(t)} - \bG_i^t ) \right\|_2^2
		\nonumber
		\\
		&\quad + \mu \E \| \bw^t - \bw^* \|_2^2
		\nonumber
		\\
		&\geq \frac{\mu}{2} \E \| \bw^t - \bw^* \|_2^2 - \frac{L^2}{2\mu n} \sum_{i=1}^{n} \E \| \bw^t - \bw^{\tau_{i}(t)} \|_2^2 
		\nonumber
		\\
		&\geq \frac{\mu}{2} \err_t - \frac{L^2T^2}{2 \mu} \eta_{t-T}^2 \left( \frac{2 \sigma^2}{n} + \ConstLemOne \err_{t}^{\max} \right)
		\nonumber
		\\
		&= - \frac{L^2 T^2}{\mu} \eta_{t-T}^2 \frac{\sigma^2}{n} - \frac{\ConstLemOne L^2 T^2}{2 \mu} \eta_{t-T}^2 \err_{t}^{\max} + \frac{\mu}{2} \err_t,
		\label{eq:blacktri}
	\end{align}
	where the last inequality also follows from Lemma \ref{lem:diff}. 
	Thus, plugging \eqref{eq:bigstar} and \eqref{eq:blacktri} back into \eqref{eq:star+tri} gives
	\begin{align}
		B \geq& 
		\left( \frac{\mu}{2} \hspace{-0.5mm}-\hspace{-0.5mm} \frac{L^2}{2\alpha} \right) \err_t - \left( \alpha \hspace{-0.5mm}+\hspace{-0.5mm} \frac{L^2}{2\alpha} \hspace{-0.5mm}+\hspace{-0.5mm} \frac{L^2}{2 \mu} \right) \ConstLemOne T^2 \eta_{t-T}^2 \err_{t}^{\max}
		\nonumber
		\\
		& -\left( \alpha + \frac{L^2}{2\alpha} + \frac{L^2}{2\mu} \right) T^2 \eta_{t-T}^2 \frac{2 \sigma^2}{n}.
		\label{eq:Bgeq}
	\end{align}
	Setting $\alpha = 2L^2 / \mu$ in \eqref{eq:Bgeq} 
	gives
	\begin{align}
		B &\geq \frac{\mu}{4} \err_t - \left( \frac{\mu}{4} + \frac{5 L^2}{2 \mu} \right) \ConstLemOne T^2 \eta_{t-T}^2 \err_{t}^{\max}
		\nonumber
		\\
		&\quad - \left( \frac{\mu}{2} + \frac{5 L^2}{\mu} \right) T^2 \eta_{t-T}^2 \frac{\sigma^2}{n}.
		\label{eq:spade}
	\end{align}
	
	Finally, plugging \eqref{eq:club} and \eqref{eq:spade} back into \eqref{eq:club+spade} gives the desired bound.
	
	\subsection{Proof of Lemma \ref{lem:diff}} \label{ap:diff}
	The proof is elementary:
	\begin{align}
		& \E \| \bw^t - \bw^{\tau_{i}(t)} \|_2^2
		= \E \left\| \sum_{s = \tau_{i}(t)}^{t-1} \left( \bw^{s+1} - \bw^{s} \right) \right\|_2^2
		\nonumber
		\\ 
		& \leq \left| t - \tau_{i}(t) \right| \sum_{s = \tau_{i}(t)}^{t-1} \E \| \bw^{s+1} - \bw^{s} \|_2^2
		\nonumber
		\\
		&\leq T \sum_{s = \tau_{i}(t)}^{t-1} \E \| (\eta_s/n) \bg^s \|_2^2
		\nonumber
		\\
		&\leq T \sum_{s = \tau_{i}(t)}^{t-1} \eta_{t-T}^2 \left( \frac{2 \sigma^2}{n} + \ConstLemOne  \max_{ s \in [ (t-T)_+, t ] } \err_s \right)
		\label{eq:xxx}
		\\
		&\leq T^2 \eta_{t-T}^2 \left( \frac{2 \sigma^2}{n} + \ConstLemOne \errMax_{t} \right).
		\nonumber
	\end{align}
	where \eqref{eq:xxx} follows from the monotonicity of $\{\eta_t\}_{t\geq0}$ and Lemma \ref{lem:aa}.

	\subsection{Proof of Theorem \ref{thm}} \label{ap:thm}
	Let $\rho_t \coloneqq  1 + 2 T + \left( \frac{\mu}{2} + \frac{5L^2}{\mu} \right) T^2 \eta_{t-T}$. Since $\gamma \geq 2T$, we have $\eta_{t-T} = \beta/(\gamma+t-T) \leq \beta/T$, which further implies that
	\[
	\rho_t \leq 1 + 2 T + \left( \frac{\mu}{2} + \frac{5L^2}{\mu} \right) \beta T = \bar{\rho}(T).
	\]
	Thus, the recurrence relation \eqref{eq:ana-recur} implies that
	\begin{align} 
		\err_{t+1} 
		\leq &\left(1 - \frac{\mu}{2} \eta_t \right) \err_t + \bar{\rho}(T) \ConstLemOne \eta_{t-T}^2 \err_{t}^{\max} 
		\nonumber
		\\
		& + \bar{\rho}(T) \eta_{t-T}^2 \frac{2 \sigma^2}{n}.
		\label{eq:recur}
	\end{align}
	
	Then, we prove \eqref{eq:rate} by induction. 
	
	\noindent
	\underline{i) Base case:}
	Since $\delta_2 =\gamma^2\err_0$, we have
	\begin{align*}
		&\err_{0} \leq \frac{\delta_1}{\gamma}\frac{\sigma^2}{n} + \frac{\delta_2}{\gamma^2}.
	\end{align*}
	\noindent
	\underline{ii) Induction step:}
	Suppose that for some $t \geq 0$, it holds that 
	\begin{align}
		&\err_s \leq \frac{\delta_1}{\gamma+s} \frac{\sigma^2}{n} + \frac{\delta_2}{(\gamma+t)^2}, ~s = 0, \dots, t,
		\nonumber
	\end{align}
	which implies that
	\begin{align}
		\err_{t}^{\max} 
		\hspace{-0.5mm}&=\hspace{-0.5mm} \max_{s\in[ (t-2T)_+,t]} \err_s
		\leq \frac{\delta_1}{\gamma+t-2T} \frac{\sigma^2}{n} + \frac{\delta_2}{(\gamma+t-2T)^2}.
		\nonumber
	\end{align}
	Together with \eqref{eq:recur}, this implies that
	\begin{align}
		&\err_{t+1} 
		\leq \left( 1 - \frac{\mu\beta/2}{\gamma+t} \right) \frac{\delta_1}{\gamma+t} \frac{\sigma^2}{n}
		\nonumber
		\\
		& \hspace{-0.5mm}+ \frac{\bar{\rho}(T)\ConstLemOne\beta^2\delta_1}{(\gamma+t-T)^2 (\gamma+t-2T)} \frac{\sigma^2}{n} + \frac{2\bar{\rho}(T)\beta^2}{(\gamma+t-T)^2} \frac{\sigma^2}{n}
		\nonumber
		\\
		&\hspace{-0.5mm}+ \hspace{-1mm}\left(\hspace{-1mm} 1 \hspace{-0.5mm}-\hspace{-0.5mm} \frac{\mu\beta/2}{\gamma+t} \hspace{-0.5mm}\right) \hspace{-0.5mm}\frac{\delta_2}{(\gamma+t)^2} \hspace{-0.5mm}+\hspace{-0.5mm} \frac{\ConstLemOne\beta^2\bar{\rho}(T)\delta_2}{(\gamma+t-T)^2 (\gamma\hspace{-0.5mm}+\hspace{-0.5mm}t\hspace{-0.5mm}-\hspace{-0.5mm}2T)^2}.
		\label{eq:speedup}
	\end{align}
	We note that
	\begin{align}
		\left( 1 - \frac{\mu\beta/2}{\gamma+t} \right) \frac{1}{\gamma+t}
		&=\frac{\gamma+t-1}{(\gamma+t)^2} - \frac{\mu\beta/2-1}{(\gamma+t)^2}
		\nonumber
		\\
		&\leq \frac{1}{\gamma+t+1} - \frac{\mu\beta/2-1}{(\gamma+t)^2},
		\label{eq:11}
	\end{align}
	where the inequality holds because
	\begin{align*}
		\frac{\gamma+t-1}{(\gamma+t)^2} \leq \frac{1}{\gamma+t+1}
		\Leftrightarrow (\gamma+t)^2 \geq (\gamma+t)^2-1.
	\end{align*}
	We also note that
	\begin{align}
		\left( 1 - \frac{\mu\beta/2}{\gamma+t} \right) \frac{1}{(\gamma+t)^2} 
		&=\frac{\gamma+t-2}{(\gamma+t)^3} - \frac{\mu\beta/2-2}{(\gamma+t)^3} 
		\nonumber
		\\
		&\leq \frac{1}{(\gamma+t+1)^2} - \frac{\mu\beta/2-2}{(\gamma+t)^3},
		\label{eq:22}
	\end{align}
	where the inequality holds because
	\begin{align*}
		\frac{\gamma+t-2}{(\gamma+t)^3} \leq \frac{1}{(\gamma+t+1)^2}
		\Leftrightarrow -3 (\gamma+t) \leq 2.
	\end{align*}
	Thus, Plugging \eqref{eq:11} and \eqref{eq:22} into \eqref{eq:speedup} gives
	\begin{align}
		&\err_{t+1} 
		\leq \frac{\sigma^2}{n} \left[ \frac{\delta_1}{\gamma+t+1} - \frac{(\mu\beta/2-1)\delta_1}{(\gamma+t)^2} \right.
		\nonumber
		\\
		&\left. + \frac{\ConstLemOne\beta^2\bar{\rho}(T)\delta_1}{(\gamma+t-T)^2 (\gamma+t-2T)} + \frac{2\beta^2\bar{\rho}(T)}{(\gamma+t-T)^2} \right]
		\nonumber
		\\
		&+\hspace{-0.7mm} \frac{\delta_2}{(\hspace{-0.5mm}\gamma\hspace{-0.7mm}+\hspace{-0.7mm}t\hspace{-0.7mm}+\hspace{-0.7mm}1\hspace{-0.5mm})^2}
		\hspace{-0.7mm}-\hspace{-0.7mm} \frac{(\hspace{-0.5mm}\mu\beta/2\hspace{-0.7mm}-\hspace{-0.7mm}2\hspace{-0.5mm})\delta_2}{(\gamma+t)^2} 
		\hspace{-0.7mm}+\hspace{-0.7mm} \frac{\ConstLemOne\beta^2\bar{\rho}(T)\delta_2}{(\hspace{-0.5mm}\gamma\hspace{-0.7mm}+\hspace{-0.7mm}t\hspace{-0.7mm}-\hspace{-0.7mm}T\hspace{-0.3mm})^2 (\hspace{-0.5mm}\gamma\hspace{-0.7mm}+\hspace{-0.7mm}t\hspace{-0.7mm}-\hspace{-0.7mm}2T\hspace{-0.3mm})^2}.
		\label{eq:Eup}
	\end{align}
	Moreover, $\gamma\geq2T$ implies that
	\begin{align}
		\frac{1}{(\gamma+t-T)^2} \leq \frac{4}{(\gamma+t)^2}.
		\label{eq:a1}
	\end{align}
	Since $\beta > 2/\mu$ and 
	$\gamma \geq 2T + \frac{16C_L\beta^2}{\mu\beta-2}\bar{\rho}(T)$,
	we have
	\begin{align}
		\frac{\ConstLemOne\beta^2\bar{\rho}(T)\delta_1}{\gamma+t-2T} \leq \left(\frac{\mu\beta}{2}-1\right) \frac{\delta_1}{8}.
		\label{eq:a2}
	\end{align}
	Besides, it follows from $\delta_1 = 32\beta^2\bar{\rho}(T)/(\mu\beta-2)$ that
	\begin{align}
		2\beta^2\bar{\rho}(T) = \left(\frac{\mu\beta}{2}-1\right) \frac{\delta_1}{8}.
		\label{eq:a3}
	\end{align}
	Combing \eqref{eq:a1}, \eqref{eq:a2}, and \eqref{eq:a3} gives
	\begin{align}
		&\frac{\ConstLemOne\beta^2\bar{\rho}(T)\delta_1}{(\gamma\hspace{-0.7mm}+\hspace{-0.7mm}t\hspace{-0.7mm}-\hspace{-0.7mm}T)^2 (\gamma\hspace{-0.7mm}+\hspace{-0.7mm}t\hspace{-0.7mm}-\hspace{-0.7mm}2T)} \hspace{-0.7mm}+\hspace{-0.7mm} \frac{2\beta^2\bar{\rho}(T)}{(\gamma\hspace{-0.7mm}+\hspace{-0.7mm}t\hspace{-0.7mm}-\hspace{-0.7mm}T)^2} 
		\hspace{-0.7mm}\leq\hspace{-0.7mm} \frac{(\mu\beta/2\hspace{-0.7mm}-\hspace{-0.7mm}1)\delta_1}{(\gamma+t)^2}.
		\label{eq:add1}
	\end{align}
	Since $\beta > 4/\mu$ and 
	$\gamma \geq 2T + \sqrt{\frac{8C_L\beta^2}{\mu\beta-4}\bar{\rho}(T)}$,
	we have
	\begin{align}
		\frac{\ConstLemOne\beta^2\bar{\rho}(T)}{(\gamma+t-2T)^2} \leq \frac{1}{4} \left(\frac{\mu\beta}{2}-2\right).
		\label{eq:a4}
	\end{align}
	Combining \eqref{eq:a1} and \eqref{eq:a4} gives
	\begin{align}
		&\frac{\ConstLemOne\beta^2\bar{\rho}(T)\delta_2}{(\gamma+t-T)^2 (\gamma+t-2T)^2} \leq \frac{(\mu\beta/2-2)\delta_2}{(\gamma+t)^2}.
		\label{eq:add2}
	\end{align}
	Plugging \eqref{eq:add1} and \eqref{eq:add2} back into \eqref{eq:Eup} yields
	\begin{align}
		\err_{t+1} \leq \frac{\delta_1}{\gamma+t+1} \frac{\sigma^2}{n} + \frac{\delta_2}{(\gamma+t+1)^2},
	\end{align}
	which completes the proof of the induction step.
	
\end{document}